\documentclass[11pt]{article}

\usepackage{latexsym,mathrsfs}
\usepackage{amsmath,amssymb}
\usepackage{amsthm,enumerate,verbatim}
\usepackage{amsfonts}
\usepackage{graphicx}
\usepackage{algorithm}
\usepackage{algorithmic}
\usepackage{url}

\newtheorem{theorem}{Theorem}

\newtheorem{corollary}{Corollary}

\newtheorem*{definition*}{Definition}

\newtheorem{remark}{Remark}

\numberwithin{equation}{section}
\numberwithin{table}{section}
\numberwithin{figure}{section}
\def \R{{\mathbb R}}
\def \C{{\mathbb C}}

\newcommand {\mat}  [1] {\left[\begin{array}{#1}}
\newcommand {\rix}      {\end{array}\right]}

\DeclareMathOperator{\argmin}{argmin}

\DeclareMathOperator{\rank}{rank}

\newcommand{\eproof}{\space
    {\ \vbox{\hrule\hbox{\vrule height1.3ex\hskip0.8ex\vrule}\hrule}}\par}

\title{A note on approximating the nearest stable discrete-time descriptor system with fixed rank}
\date{}

\author{Nicolas Gillis\thanks{Department of Mathematics and Operational Research,
Facult\'e Polytechnique, Universit\'e de Mons, Rue de Houdain~9, 7000 Mons, Belgium; \texttt{nicolas.gillis@umons.ac.be}. N. Gillis
acknowledges the support of the ERC (starting grant n$^\text{o}$ 679515) and F.R.S.-FNRS (incentive grant for scientific research n$^\text{o}$ F.4501.16).}
 \qquad Michael Karow\thanks{TU Berlin, Institut f${\rm \ddot{u}}$r Mathematik,
Stra{\ss}e des 17. Juni 136, 10623 Berlin, Germany; \texttt{karow@math.tu-berlin.de. }} 
 \qquad Punit Sharma\thanks{Department of Mathematics, Indian Institute of Technology Delhi, Hauz Khas, New Delhi-110016, India; \texttt{punit.sharma@maths.iitd.ac.in}.} 
}

\begin{document}

\maketitle

\begin{abstract}
Consider a discrete-time linear time-invariant descriptor system  $Ex(k+1)=Ax(k)$ for $k \in \mathbb Z_{+}$. In this paper, we tackle for the first time the problem of stabilizing such systems by computing a nearby regular index one stable system $\hat E x(k+1)= \hat A x(k)$ with $\text{rank}(\hat E)=r$. We reformulate this highly nonconvex problem into an equivalent optimization problem with a relatively simple feasible set onto which it is easy to project. This allows us to employ a block coordinate descent method to obtain a nearby regular index one stable system. We illustrate the effectiveness of the algorithm on several examples.
\end{abstract}


\section{Introduction}

In~\cite{OrbNV13, GilKS18a}, authors have tackled the problem of computing the nearest stable matrix in the discrete case, that is, given an unstable matrix $A$, find the smallest perturbation $\Delta_A$ with respect to Frobenius norm such that
$\hat A = A+\Delta_A$ has all its eigenvalues inside the unit ball centred at the origin. In this paper, we aim to
generalize the results in~\cite{GilKS18a} for matrix pairs $(E,A)$, where $E,A\in \R^{n,n}$.
%

The matrix pair $(E,A)$ is called \emph{regular} if
$\operatorname{det}(\lambda E-A)\neq 0$ for some $\lambda \in \mathbb C$, which we denote $\operatorname{det}(\lambda E-A) \not\equiv 0$, otherwise it is called \emph{singular}.
For a regular matrix pair $(E,A)$, the roots of the polynomial $\operatorname{det}(z E-A)$ are called \emph{finite eigenvalues} of the pencil
$zE-A$ or of the pair $(E,A)$. A regular pair $(E,A)$ has \emph{$\infty$ as an eigenvalue} if $E$ is singular.
A regular real matrix pair $(E,A)$ 
can be transformed to \emph{Weierstra\ss\ canonical form} \cite{Gan59a}, that is, there exist nonsingular matrices $W, T \in \C^{n,n}$ such that
\[
E=W\mat{cc}I_q& 0\\0&N\rix T \quad \text{and}\quad A=W \mat{cc}J &0\\0&I_{n-q}\rix T,
\]
where $J \in \C^{q,q}$ is a matrix in \emph{Jordan canonical form} associated with the $q$ finite eigenvalues of
the pencil $z E-A$ and  $N \in \C^{n-q,n-q}$ is a nilpotent matrix in Jordan canonical form  corresponding
to  $n-q$ times the  eigenvalue $\infty$. If $q < n$ and $N$ has degree of nilpotency $\nu \in \{1,2,\ldots\}$, that is, $N^{\nu}=0$ and $N^i \neq 0$ for $i=1,\ldots,\nu-1$, then $\nu$ is called the \emph{index of the pair} $(E,A)$. If $E$
is nonsingular, then by convention the index is $\nu=0$; see for example~\cite{Meh91,Var95}.
The  matrix pair $(E,A) \in (\R^{n,n})^2$ is said to be \emph{stable} (resp.\@ \emph{asymptotically stable})
if all the finite eigenvalues of $zE-A$ are
in the closed (resp.\@ open) unit ball and those on the
unit circle are semisimple.
The  matrix pair $(E,A)$ is said to be \emph{admissible} if it is regular, of index at most one, and
stable.

The various distance problems for linear  control systems is an important research topic in the numerical linear algebra community; for example,
the distance to bounded realness~\cite{AlaBKMM11},
the robust stability problem~\cite{Zho11},
the stability radius problem for standard systems~\cite{Bye88,HinP86}
and for descriptor systems~\cite{ByeN93,DuLM13},
the nearest stable matrix problem
for continuous-time systems~\cite{OrbNV13,GilS17,MehMS17,GugL17} and
for discrete-time systems~\cite{OrbNV13,NesP17,GP2018,GilKS18a},
the nearest continuous-time admissible descriptor system
problem~\cite{GilMS17}, and the nearest positive real system problem~\cite{GilS17b}.

For a given unstable matrix pair $(E,A)$, the discrete-time nearest stable matrix pair problem is to solve the
following optimization problem
\begin{equation}\label{mainprob}
\inf_{(\hat E,\hat A)\in\mathcal S^{n,n}} {\|E-\hat E\|}_F^2+{\|A-\hat A\|}_F^2,\tag{$\mathcal{P}$}
\end{equation}
where $\mathcal S^{n,n}$ is the set of admissible pairs of size $n \times n$.
This problem is the converse of stability radius problem for descriptor systems~\cite{ByeN93,DuLM13} and
the discrete-time counter part of continuous-time nearest stable matrix pair problem~\cite{GilMS17}.
Such problems arise in systems identification where one needs to identify a stable matrix pair depending on observations~\cite{OrbNV13,GilS17}.
This is a highly nonconvex optimization problem because the set $\mathcal S^{n,n}$ is unbounded, nonconvex  and neither open nor closed. In fact,
consider the matrix pair
\begin{equation}\label{eq:ex1}
(E,A)=\Bigg(
\mat{ccc}1&0&0\\0&0&0\\0&0&0 \rix,~\mat{ccc}1/2&0&2\\0&1&0\\0&0&1 \rix
\Bigg).
\end{equation}
The pair $(E,A)$ is regular since $\text{det}(\lambda E-A)=\text{det}(\lambda -1/2)\not\equiv 0$, of index one,
and stable with the only finite eigenvalue $\lambda_1=1/2$. Thus $(E,A) \in \mathcal S^{3,3}$.
Let
\begin{equation}\label{eq:ex1perturb}
(\Delta_E,\Delta_A)=\Bigg(
\mat{ccc}0&0&0\\0&\epsilon_1&\epsilon_2\\0&0&0 \rix,\mat{ccc}0&0&0\\0&0&0\\0&0&-\delta \rix
\Bigg),
\end{equation}
and consider the perturbed pair $(E+\Delta_E,A+\Delta_A)$. If we let $\delta=\epsilon_1=0$ and
$\epsilon_2>0$, then the perturbed pair is still regular and stable as the only finite eigenvalue $\lambda_1=1/2$
belongs to the unit ball, but it is of index two.
For $\epsilon_2=\delta=0$ and $0<\epsilon_1<1$, the perturbed pair is regular,
of index one but has two finite eigenvalues $\lambda_1=1/2$ and $\lambda_2=1/\epsilon_1 >1$. This implies that
the perturbed pair is unstable. This shows that $\mathcal S^{3,3}$ is not open. Similarly, if we let $\epsilon_1=\epsilon_2=0$ and
$\delta >0$, then as $\delta \rightarrow 1$ the perturbed pair becomes non-regular. This shows that $\mathcal S^{3,3}$ is not closed.
The nonconvexity of $\mathcal S^{n,n}$ follows by considering for example
\begin{equation}\label{eq:nonconvex}
\Sigma_1=\Big(I_2,\underbrace{\mat{cc}0.5 & 2\\ 0& 1 \rix}_{A}\Big), \quad \Sigma_2=\Big(I_2,\underbrace{\mat{cc}0.5 & 0\\ -2& 1\rix}_{B}\Big),
\end{equation}
where $\Sigma_1,\Sigma_2 \in \mathcal S^{2,2}$, while $\gamma \Sigma_1 + (1-\gamma)\Sigma_2 \notin \mathcal S^{2,2}$
for $\gamma=\frac{1}{2}$, since $\frac{1}{2} \Sigma_1+\frac{1}{2} \Sigma_2$ has two eigenvalues 0.75$\pm$0.96$i$ outside the unit ball.
Therefore it is in general difficult to work directly with the set $\mathcal S^{n,n}$.
We explain in Section~\ref{reform} the difficulty in generalizing the results in~\cite{GilKS18a} for problem~\eqref{mainprob}.

In this paper, we consider instead a \emph{rank-constrained nearest stable matrix pair problem}. For this, let $r (<n) \in \mathbb Z_{+}$ and let us
define a subset $\mathcal S_r^{n,n}$ of $\mathcal S^{n,n}$ by
\[
\mathcal S_r^{n,n} :=\left\{(\hat E,\hat X)\in \mathcal S^{n,n}:~\text{rank}(\hat E)=r\right\}.
\]
For a given unstable matrix pair $(E,A)$, the rank-constrained nearest stable matrix pair problem requires to compute the smallest perturbation $(\Delta_E,\Delta_A)$
with respect to Frobenius norm such that $(E+\Delta_E,A+\Delta_A)$ is admissible with $\text{rank}(E+\Delta_E)=r$, or equivalently, we aim to solve the following
optimization problem
\begin{equation}\label{restprob}
\inf_{(\hat E, \hat A)\in\mathcal S_r^{n,n}} {\|E-\hat E\|}_F^2+{\|A-\hat A\|}_F^2 \tag{$\mathcal{P}_r$}.
\end{equation}
The problem~\eqref{restprob} is also nonconvex as the set $\mathcal S_r^{n,n}$ is nonconvex.
To solve~\eqref{restprob}, we provide a simple parametrization of $\mathcal S_r^{n,n}$
in terms of a matrix quadruple $(T,W,U,B)$, where $T,W\in \R^{n,n}$ are invertible, $U\in \R^{r,r}$ is
orthogonal, and $B\in \R^{r,r}$ is a positive semidefinite contraction, see Section~\ref{reform}.
This parametrization
results in an equivalent optimization problem with a feasible set onto which it is easy to project, and we derive a block coordinate descent method to tackle it; see Section~\ref{sec:algo}. We illustrate the effectiveness of our algorithm
over several numerical examples in Section~\ref{sec:numexp}.

\paragraph{Notation}
Throughout the paper, $X^T$ and
$\|X\|$  stand for the transpose and the spectral norm of a real square matrix $X$, respectively.
We write $X\succ 0$ and $X\succeq 0$ $(X \preceq 0)$ if $X$ is symmetric and positive definite
or positive semidefinite (symmetric negative semidefinite), respectively.
By $I_m$ we denote the identity matrix of size $m \times m$.

\section{Reformulation of problem~\eqref{restprob}} \label{reform}

As mentioned earlier the set $\mathcal S_r^{n,n}$ is nonconvex. It is also an unbounded set which is
 neither open nor closed.
 Consider
\begin{equation*}
\tilde\Sigma_1=\Big(\mat{cc}I_r & 0\\ 0& 0 \rix,\mat{cc}A_1 & 0\\ 0& I_{n-r} \rix\Big), \quad
\tilde\Sigma_2=\Big(\mat{cc}I_r & 0\\ 0& 0 \rix,\mat{cc}B_1 & 0\\ 0& I_{n-r} \rix\Big),
\end{equation*}
where $A_1=\mat{cc}A &0\\0& I_{r-2}\rix$, $B_1=\mat{cc}B &0\\0& I_{r-2}\rix$, and $A$ and $B$ are defined as in~\eqref{eq:nonconvex}.
We have that $\tilde\Sigma_1,\tilde\Sigma_2 \in \mathcal S_r^{n,n}$ because $(I_r,A_1)$ and $(I_r,B_1)$ are stable.
Moreover $\frac{1}{2} \tilde\Sigma_1+\frac{1}{2} \tilde\Sigma_2 \notin \mathcal S_r^{n,n}$  as it has two eigenvalues 0.75$\pm$0.96$i$ outside the unit ball hence $\mathcal S_r^{n,n}$ is non-convex. 
 To show that $\mathcal S_r^{n,n}$ is neither open nor
closed, let $(E,A)$ and $(\Delta_E,\Delta_A)$ be as defined in~\eqref{eq:ex1} and~\eqref{eq:ex1perturb}, and consider
\[
(\tilde E,\tilde A)=\Big(\mat{cc}I_{r-1} & 0\\ 0& E \rix,\mat{cc}I_{r-1} & 0\\ 0& A \rix\Big)
\]
and the perturbation
\[
(\Delta_{\tilde E},\Delta_{\tilde A})=\Big(\mat{cc}I_{r-1} & 0\\ 0& \Delta_E \rix,\mat{cc}I_{r-1} & 0\\ 0& \Delta_A \rix\Big).
\]
By using similar arguments as in the case of $\mathcal S^{n,n}$ one can show that $\mathcal S_r^{n,n}$
is neither open nor closed. Therefore it is difficult to compute a global solution to problem~\eqref{restprob}
and to work directly with the set $\mathcal S_r^{n,n}$. For this reason, we reformulate the rank-constrained nearest stable matrix pair problem into an equivalent problem with a relatively simple feasible set.
In order to do this , we derive a parametrization of admissible pairs into invertible, symmetric and orthogonal matrices.
We first recall a result from~\cite{GilKS18a} that gives a characterization for stable matrices.
\begin{theorem}{\rm \cite[Theorem 1]{GilKS18a}}\label{thm:stabmatchar}
Let $A\in \R^{n,n}$. Then $A$ is stable if and only if $A=S^{-1}UBS$ for some $S,U,B \in \R^{n,n}$
such that $S\succ 0$, $U^TU=I_n$, $B\succeq 0$, and $\|B\|\leq 1$.
\end{theorem}
We note that, in the proof of Theorem~\ref{thm:stabmatchar}, only the invertibility of matrix
$S$ is needed and the condition of symmetry on $S$ can be relaxed. We found that
this relaxation on matrix $S$ does not make any difference on the numerical results in~\cite{GilKS18a}.
The only gain is that the projection of $S$ on the set of positive definite matrices takes some
time and that can be avoided. Therefore, we rephrase the definition of a SUB matrix in~\cite{GilKS18a}
and the corresponding characterization of stable matrices as follows.
\begin{theorem}\label{thm:newstabmatchar}
Let $A\in \R^{n,n}$. Then $A$ is stable if and only if $A$ admits a SUB form, that is,
 $A=S^{-1}UBS$ for some $S,U,B \in \R^{n,n}$ such that $S$ is invertible, $U^TU=I_n$, $B\succeq 0$, and $\|B\|\leq 1$.
\end{theorem}
\begin{theorem}\label{thm:reform1}
Let $E,A\in \R^{n,n}$ be such that $\text{rank}(E)=r$. Then
$(E,A)$ is admissible if and only if
there exist matrices $T,W \in \R^{n,n}$, $S,U,B\in \R^{r,r}$ such that the matrices $T,W,S$ are invertible,
$U^TU=I_r$, $B\succeq 0$, $\|B\|\leq 1$ such that
\begin{equation}\label{eq:firstreform}
E=W \begin{bmatrix}   I_r  & 0 \\   0 & 0 \end{bmatrix}T, \quad
\text{and}
\quad
A=W\begin{bmatrix}  S^{-1}UBS  & 0 \\  0 & I_{n-r}\end{bmatrix}T.
\end{equation}
\end{theorem}
\proof For a regular index one pair $(E,A)$,
there exist invertible matrices $W,T\in \R^{n,n}$ such that
\begin{equation}
E=W \begin{bmatrix}   I_r  & 0 \\   0 & 0 \end{bmatrix}T \qquad \text{ and } \qquad
A=W\begin{bmatrix}  \tilde A  & 0 \\  0 & I_{n-r}\end{bmatrix}T,
\end{equation}
see~\cite{Dai89}.
Further, the finite eigenvalues of $(E,A)$ and $\tilde A$ are same because $\text{det}(\lambda E-A)=0$
if and only if
 $\text{det}(\lambda I_r-\tilde A)=0$.
Thus by stability of $(E,A)$ and Theorem~\ref{thm:newstabmatchar}, it follows that $\tilde A$ admits a SUB form, that is,
there exist $S,U,B \in \R^{r,r}$ such that
$S$ is invertible, $U^TU=I_r$, $B\succeq 0$, $\|B\|\leq 1$, and $\tilde A =S^{-1}UBS$.\\
Conversely, it is easy to see that any matrix pair $(E,A)$ in the form~\eqref{eq:firstreform} is regular and of index one.
The stability of $(E,A)$ follows from Theorem~\ref{thm:newstabmatchar} as the matrix $S^{-1}UBS$ is stable.
\eproof
If the matrix $E$ is nonsingular, then Theorem~\ref{thm:reform1}
can be further simplified as follows.
\begin{theorem}\label{eq:nonsingcase}
Let $E,A \in \R^{n,n}$, and let $E$ be nonsingular. Then $(E,A)$
is admissible if and only if there exist
matrices  $S,U,B\in \R^{n,n}$ such that $A=S^{-1}UBSE$, where
 $S$ is invertible, $U^TU=I_n$, $B\succeq 0$, and $\|B\|\leq 1$.
\end{theorem}
\proof
Since $E$ is nonsingular, the matrix pair $(E,A)$ can be equivalently written as a standard pair $(I_n,AE^{-1})$,
and then stability of $(E,A)$ can be determined by the eigenvalues of $AE^{-1}$. That means, $(E,A)$ is stable
if and only if $AE^{-1}$ is stable. Thus from Theorem~\ref{thm:newstabmatchar}, $AE^{-1}$ is stable if and only if
$AE^{-1}$ admits a SUB form, that is, $AE^{-1}=S^{-1}UBS$ for some $S,U,B \in \R^{n,n}$ such that $S$ is invertible,
$U^TU=I_n$, $B\succeq 0$ and $\|B\|\leq 1$.
\eproof
We note that, for a standard pair $(I_n,A)$ (with $E=I_n$), Theorem~\ref{thm:reform1} coincides with
Theorem~\ref{thm:newstabmatchar} as in this case $W$ and $T$ can be chosen to be the identity matrix which yields $A=S^{-1}UBS$.
A similar result also holds for asymptotically stable matrix pairs which can be seen as a generalization of \cite[Theorem 2]{GilKS18a}.
\begin{theorem}
Let $E,A\in \R^{n,n}$ be such that $\text{rank}(E)=r$. Then
$(E,A)$ is regular, of index one and asymptotically stable if and only if
there exist matrices $T,W \in \R^{n,n}$, $S,U,B\in \R^{r,r}$ such that the matrices $T,W,S$ are invertible,
$U^TU=I_r$, $B\succeq 0$, $\|B\|< 1$ such that
\begin{equation} \label{eq:secondreform}
E=W \begin{bmatrix}   I_r  & 0 \\   0 & 0 \end{bmatrix}T, \quad
\text{and}
\quad
A=W\begin{bmatrix}  S^{-1}UBS  & 0 \\  0 & I_{n-r}\end{bmatrix}T.
\end{equation}
\end{theorem}
\proof The proof follows is similar to that of Theorem~\ref{thm:reform1} by using \cite[Theorem 2]{GilKS18a} instead
of Theorem~\ref{thm:newstabmatchar}.
\eproof
Note that the matrix $S$ is invertible in Theorem~\ref{thm:reform1} and therefore it can be absorbed in $W$ and $T$. The advantage is that this reduces the number of variables in the corresponding optimization problem.
%
\begin{corollary}\label{thm:reform2}
Let $E,A\in \R^{n,n}$ be such that $\text{rank}(E)=r$. Then
$(E,A)$ is admissible if and only if
there exist invertible matrices $T,W \in \R^{n,n}$, and $U,B\in \R^{r,r}$ with
$U^TU=I_r$, $B\succeq 0$ and $\|B\|\leq 1$ such that
\begin{equation}\label{eq:secreform}
E=W \begin{bmatrix}   I_r  & 0 \\   0 & 0 \end{bmatrix}T, \quad
\text{and}
\quad
A=W\begin{bmatrix}  UB & 0 \\  0 & I_{n-r}\end{bmatrix}T.
\end{equation}
\end{corollary}
In view of Corollary~\ref{thm:reform2}, the set $\mathcal S_r^{n,n}$ of restricted rank admissible pairs
can be characterized in terms of matrix pairs~\eqref{eq:secreform}, that is,
\begin{eqnarray*}
&\mathcal S_r^{n,n} = \Bigg \{\left(W \mat{cc}  I_r  & 0 \\   0 & 0 \rix T,
W\mat{cc}  UB  & 0 \\  0 & I_{n-r}\rix T\right):~\text{invertible}~ T,W \in \R^{n,n},\\
& \hspace{6cm} U,B\in\R^{r,r},U^TU=I_r,B\succeq 0,\|B\|\leq 1 \Bigg \}.
\end{eqnarray*}
This parametrization changes the feasible set and the objective function in problem~\eqref{restprob} as
\begin{equation}\label{eq:reform_prob}
(\mathcal{P}_r) \quad = \quad  \inf_{W,T \in \R^{n,n},\,UB\in \R^{r,r},\,U^TU=I_r,\, \|B\|\leq 1} \;  f(W,T,U,B),
\end{equation}
where
\[
f(W,T,U,B)=
{\left\|E-W \mat{cc}  I_r  & 0 \\   0 & 0 \rix T\right\|}_F^2+
{\left\|A-W\mat{cc}  UB  & 0 \\  0 & I_{n-r}\rix T\right\|}_F^2 .
\]
An advantage of this reformulation over~\eqref{restprob} is that it is relatively easy to project onto the feasible set of~\eqref{eq:reform_prob}.
This enables us to use standard optimization schemes to solve it, see Section~\ref{sec:algo}.

As mentioned in~\cite{GilMS17}, for the standard pair $(I_n,A)$ making $A$ stable without perturbing the identity matrix gives an upper
bound to the solution of $({\mathcal P}_n)$, because
\begin{equation}\label{eq:numupbound}
\inf_{(M,X)\in \mathcal S_n^{n,n}} {\|I_n-M\|}_F^2+{\|A-X\|}_F^2\leq \inf_{(I_n,X)\in \mathcal S_n^{n,n}} {\|A-X\|}_F^2
=\inf_{(I_n,S^{-1}UBS)\in \mathcal S_n^{n,n}} {\|A-S^{-1}UBS\|}_F^2.
\end{equation}
Note that the right hand side infimum in~\eqref{eq:numupbound} is the distance of $A$ from the set of stable matrices~\cite{GilKS18a}.
It is demonstrated in our numerical experiments that (as expected) the inequality in~\eqref{eq:numupbound} is strict. We also note that
similar arguments do not extend to the solution of problem~\eqref{restprob}, when $r < n$. In this case, the distance of $A$ from the set of stable matrices
is not an upper bound for the solution of~\eqref{restprob}, see Section~\ref{sec:numexp}.
We close the section with a remark that emphasizes the difficulty in solving~\eqref{mainprob} over~\eqref{restprob}.

\begin{remark}{\rm
In view of Corollary~\ref{thm:reform2}, the set $\mathcal S^{n,n}$ of admissible pairs can be written as
\begin{eqnarray*}
\mathcal S^{n,n}
&=&\bigcup_{r=1}^n \mathcal S_r^{n,n}.
\end{eqnarray*}
Hence we have that
\[
\eqref{mainprob} = \min_{r=1,2,\ldots,n}\eqref{restprob}.
\]
To compute a solution of~\eqref{mainprob},
a possible way is therefore to solve $n$ rank-constrained problems~\eqref{restprob}.
For $n$ large, this would be rather costly as it makes the corresponding algorithm for~\eqref{mainprob} $n$ times more expensive than for~\eqref{restprob}.
However, in practice, the rank $r$ has to be chosen close to the (numerical) rank of $E$ so that it can be estimated from the input data.
Also, as we will see in Section~\ref{sec:numexp}, the error tends to change monotonically with $r$ (first it decreases as $r$ increases --unless $r=1$ is the best value-- and then increases after having achieved the best value for $r$) which could also be used to avoid computing the solutions for all $r$.
}
\end{remark}

\section{Algorithmic solution for~\eqref{restprob}} \label{sec:algo}

To solve~\eqref{eq:reform_prob}, we use a block coordinate descent method and optimize alternatively over $W$, $T$ and $(U,B)$. For $T$, $U$ and $B$ fixed, the optimal $W$ can be computed using least squares, and similarly for the optimal $T$. Note that the least squares problem in $W$ (resp.\@ $T$) can be solved independently for each row (resp.\@ each column)
To update $(U,B)$ for $W$ and $T$ fixed, we use the fast gradient method from~\cite{GilKS18a} (it can be easily adapted by fixing $S$ to the identity and modifying the gradients).

\algsetup{indent=2em}
\begin{algorithm}[ht!]
\caption{Block Coordinate Descent Method for~\eqref{eq:reform_prob}} \label{bcd}
\begin{algorithmic}[1]
\REQUIRE
An initialization
$W \in \mathbb{R}^{n \times n}, T \in  \mathbb{R}^{n \times n}, U \in \mathbb{R}^{r \times r}, B \in \mathbb{R}^{r \times r}$.

\ENSURE An approximate solution $(W,T,U,B)$ to~\eqref{eq:reform_prob}.  \medskip

\FOR{$k = 1, 2, \dots$}

\STATE $W \leftarrow \argmin_{Y} f(Y,T,U,B)$;  \emph{\% Least squares problem}

\STATE $T \leftarrow \argmin_{X} f(W,X,U,B)$; \emph{\% Least squares problem}

\STATE Apply a few steps of the fast gradient method from~\cite{GilKS18a} on
\[
\min_{(U,B) \text{ s.t. } U^TU=I_r, \|B\| \leq 1} f(W,T,U,B)
\]
to update $(U,B)$.

\ENDFOR

\end{algorithmic}
\end{algorithm}

\subsection{Initialization}

For simplicity, we only consider one initialization scheme in this paper which is similar to the one that performed best in~\cite{GilKS18a}.
However, it is important to keep in mind that Algorithm~\ref{bcd} is sensitive to initialization and that coming up with good initialization schemes is a topic of further research.

We take $W=T=I_n$ and $(U,B)$ as the optimal solution of
\[
\min_{(U,B) \text{ s.t. } U^TU=I_r, \|B\| \leq 1} {\|A_{1:r,1:r}-UB\|}_F^2.
\]
In this particular case, it can be computed explicitly using the polar decomposition of $A_{1:r,1:r}$~\cite{GilKS18a}.

\section{Numerical experiments} \label{sec:numexp}

In this section, we apply Algorithm~\ref{bcd} on several examples. As far as we know, there does not exist any other algorithm to stabilize matrix pairs (in the discrete case) hence we cannot compare it to another technique.  However, when $E=I_n$, we will compare to the fast gradient method of~\cite{GilKS18a} which provides a nearby stable matrix (but does not allow to modify $E$).
Our code is available from \url{https://sites.google.com/site/nicolasgillis/} and the numerical examples presented below can be directly run from this online code.
All tests are preformed using Matlab
R2015a on a laptop Intel CORE i7-7500U CPU @2.7GHz 24Go RAM.
Algorithm~\ref{bcd} runs in $O(n^3)$ operations per iteration,
including projections onto the set of orthogonal matrices, the resolution of the least squares problem and all necessary matrix-matrix products.
Hence Algorithm~\ref{bcd} can be applied on a standard laptop with $n$ up to a thousand (each iteration on the specified laptop takes about 10 seconds for $r=n$).

\subsection{Grcar matrix}

\begin{figure}[ht!]
\begin{center}
\includegraphics[width=\textwidth]{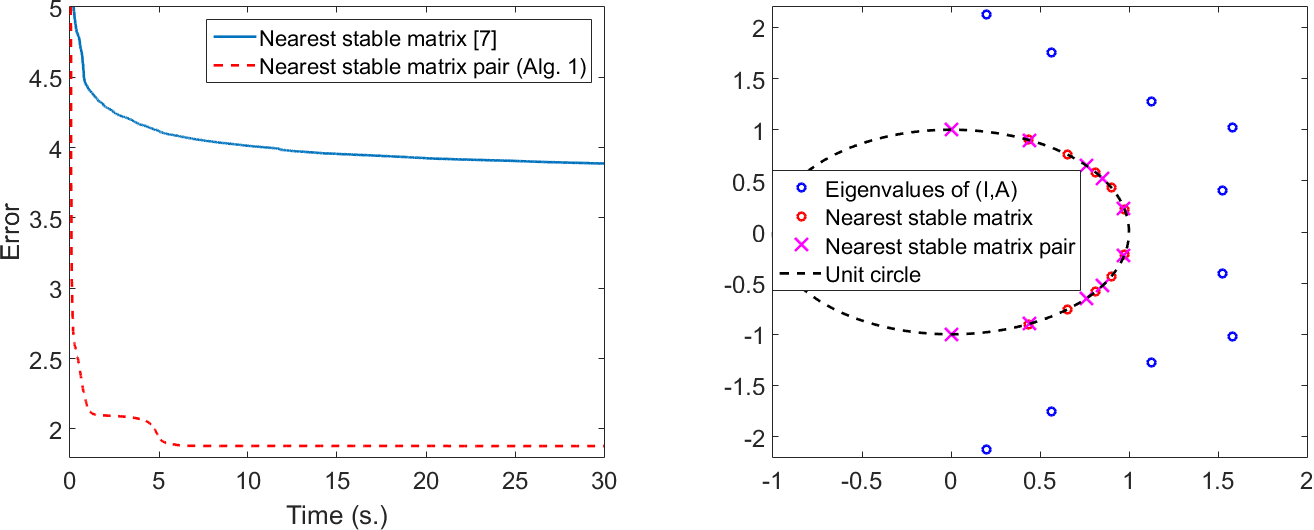}
\caption{
(Left) Evolution of the error $\|E-\hat E\|_F^2 + \|A-\hat A\|_F^2$ for the Grcar matrix of dimension 10 and order 3 in the matrix and matrix pair cases (in the matrix case, $\hat E = I_n$).
(Right) Location of the eigenvalues of $A$, and of the solutions in the matrix case and in the matrix pair case.
 \label{fig:grcar}}
\end{center}
\end{figure}

Let us first consider the pair $(I_n,A)$ where $A$ is the Grcar matrix of dimension $n$ and order $k$~\cite{GilS17}. For $n = 10$ and $k=3$,
the nearest stable matrix found in~\cite{GilKS18a} has relative error
${\|A-\hat A\|}_F^2 = 3.88$.
Applying Algorithm~\ref{bcd} with $r=n$, we obtain a matrix pair $(\hat E, \hat A)$ such that
${\|A-\hat A\|}_F^2 + {\|E-\hat E\|}_F^2  = 1.88$.
Figure~\ref{fig:grcar} displays the evolution of the error (left) and the eigenvalues of the solutions (right).
We observe that allowing $\hat E$ to be different than the identify matrix allows the matrix pair $(\hat E, \hat A)$ to be much closer to $(I_n,A)$ and have rather different eigenvalues.

\paragraph{Effect of the dimension $n$} Let us perform the same experiment as above except that we increase the value of $n$. Table~\ref{tab:grcarn} compares the error of the nearest stable matrix and of the nearest stable matrix pair. As $n$ increases, the nearest stable matrix pair allows to decrease the error of approximation.

\begin{center}
 \begin{table}[h!]
 \begin{center}
\caption{Comparison of the error for Grcar matrices $A$ of order $k=3$ and $E = I_n$ for different values of $n$.} \label{tab:grcarn}
 \begin{tabular}{|c||c|c|c|c|c|}
 \hline
                 & $n= 5$ & $n= 10$  & $n= 20$ & $n= 50$ & $n= 100$ \\
                 & (30 s.)& (60 s.)& (120 s.) & (300 s.) & (600 s.)\\ \hline
 Stable matrix   &  1.76  & 3.88   &  15.89   & 68.18    & 160.00   \\
 Stable pair     &  1.16  & 1.88   &  3.02    & 8.69     &  20.41\\ \hline
\end{tabular}
 \end{center}
 \end{table}
 \end{center}

\paragraph{Effect of $r=\rank(\hat E)$ and $\rank(E)$}
Let us now perform more extensive numerical experiments on the Grcar matrix of dimension $n=10$ of order $3$.
Let us fix $0 \leq p \leq n-1$ and define $E(i,i) = 1$ for $i > p$ otherwise $E(i,j) = 0$ (that is, $E$ is the identity matrix where $p$ diagonal entries have been set to zero) with $\rank(E) = n-p$.
Table~\ref{tab:grcarerr} gives the error of the solution obtained by Algorithm~\ref{bcd} for $r=1,2,\dots,n$.
\begin{center}
 \begin{table}[h!]
 \begin{center}
\caption{Comparison of the error for Grcar matrices $A$ with $k=3$,
and $E(i,i) = 1$ for $i > p$ otherwise $E(i,j) = 0$.} \label{tab:grcarerr}
 \begin{tabular}{|c||cccccccccc|}
 \hline
$\rank(E)$ & $r= 1$& $r= 2$& $r= 3$& $r= 4$& $r= 5$& $r= 6$& $r= 7$& $r= 8$& $r= 9$& $r=10$ \\ \hline \hline
   10      & 9.02  & 8.05  & 7.09  & 6.20  & 5.44  & 4.63  & 3.94  & 3.16  & 2.16  & \textbf{1.88}  \\ \hline
    9      & 8.04  & 7.08  & 6.13  & 5.33  & 4.61  & 3.83  & 3.16  & 2.16  & 1.57  & \textbf{1.36}  \\ \hline
    8      & 7.05  & 6.10  & 5.17  & 5.16  & 4.16  & 3.16  & 2.16  & 1.46  & \textbf{1.37}  & 1.42  \\ \hline
    7      & 6.05  & 5.13  & 4.22  & 4.16  & 2.83  & 2.17  & 1.57  & 1.52  & \textbf{1.44}  & 1.44  \\ \hline
    6      & 5.07  & 4.17  & 3.27  & 3.16  & 2.04  & \textbf{1.32}  & 1.34  & 1.69  & 1.69  & 1.91  \\ \hline
    5      & 4.09  & 3.24  & 2.34  & 1.76  & \textbf{1.20}  & 1.52  & 1.30  & 1.56  & 1.74  & 3.13  \\ \hline
    4      & 3.12  & 2.26  & 1.49  & 1.25  & \textbf{1.19}  & 1.24  & 1.30  & 1.68  & 2.96  & 2.95  \\ \hline
    3      & 2.13  & 1.31  & \textbf{0.69}  & 1.19  & 1.23  & 1.29  & 1.69  & 2.83  & 2.83  & 2.83  \\ \hline
    2      & 1.15  & \textbf{0.41}  & 1.06  & 1.22  & 1.27  & 1.27  & 1.91  & 2.79  & 2.79  & 2.79  \\ \hline
    1      & \textbf{0.17}  & 0.81  & 1.21  & 1.21  & 1.36  & 1.22  & 2.71  & 2.72  & 2.72  & 4.33  \\ \hline
\end{tabular}
 \end{center}
 \end{table}
 \end{center}

We observe that
\begin{itemize}
\item In 6 out of the 10 cases, using $r = \rank(E)$ provides the best solution.
In 3 out of the 10 cases, using $r = \rank(E)+1$ provides the best solution, and in one case $r = \rank(E)+2$ provides the best solution. This illsutartes the fact that the best value for $r$ should be close to the (numerical) rank of $E$. (Of course, since we use a single initialization, there is no guarantee that the error in Table~\ref{tab:grcarerr} is the smallest possible.)

\item In all cases, the error behaves monotonically, that is, it increases as the value of $r$ goes away from the best value.

\end{itemize}

The two observations above could be used in practice to tune effectively the value of $r$: start from a value close to the numerical rank of $E$, then try nearby values until the error increases.

Table~\ref{tab:grcartim} gives the computational time for the different cases.
We use the following stopping criterion:
\[
e(i) - e(i+1) < 10^{-8} e(i),
\]
where $e(i)$ is the error obtained at the $i$th iteration, and a time limit of 60 seconds.
\begin{center}
 \begin{table}[h!]
 \begin{center}
\caption{ Time in seconds to compute the solution obtained in Table~\ref{tab:grcarerr}.
The time limit is 60 seconds. }
 \label{tab:grcartim}
 \begin{tabular}{|c||cccccccccc|}
 \hline
$\rank(E)$ & $r= 1$& $r= 2$& $r= 3$& $r= 4$& $r= 5$& $r= 6$& $r= 7$& $r= 8$& $r= 9$& $r=10$ \\ \hline \hline
   10      & 0.59  & 1.00  & 0.72  & 0.56  & 2.44  & 2.20  & 2.08  & 43.92  & 13.41  & 49.16  \\ \hline
    9      & 0.23  & 0.58  & 0.44  & 2.34  & 2.53  & 17.89  & 3.25  & 24.19  & 43.34  & 47.72  \\ \hline
    8      & 0.22  & 0.89  & 0.42  & 0.36  & 13.16  & 1.16  & 16.34  & 3.67  & 8.09  & 11.02  \\ \hline
    7      & 0.16  & 0.81  & 0.77  & 0.61  & 13.63  & 10.20  & 60  & 60  & 60  & 60  \\ \hline
    6      & 0.19  & 0.67  & 0.73  & 1.19  & 27.42  & 60  & 60  & 60  & 60  & 60  \\ \hline
    5      & 0.17  & 0.83  & 0.59  & 4.22  & 60  & 60  & 39.17  & 60  & 60  & 60  \\ \hline
    4      & 0.25  & 1.20  & 0.39  & 29.11  & 60  & 21.80  & 60  & 60  & 40.81  & 60  \\ \hline
    3      & 0.14  & 0.28  & 0.13  & 31.08  & 60  & 60  & 60  & 55.08  & 60  & 60  \\ \hline
    2      & 0.06  & 0.17  & 17.50  & 60  & 60  & 60  & 60  & 60  & 60  & 60  \\ \hline
    1      & 0.02  & 23.56  & 26.61  & 48.11  & 60  & 60  & 37.95  & 60  & 60  & 60  \\ \hline
\end{tabular}
 \end{center}
 \end{table}
 \end{center}

We observe that the algorithm converges much faster when $r$ is small.
This can be partly explained by the smaller number of variables, being $2n^2 + 2r^2$.

\subsection{Scaled all-one matrix}

In this section, we perform a similar experiment than in the previous section
with $A = \alpha ee^T$ where $e$ is the vector of all ones, which is an example from~\cite{GP2018}.
For $\alpha > 1/n$, the matrix is unstable.
For $1/n \leq \alpha \leq 2/n$,
the nearest stable matrix is $ee^T/n$.

Let us take $n=10$ and $\alpha = 2/n = 0.2$
for which the nearest stable matrix is $A = 0.1 ee^T$ with error 1.
The nearest stable matrix pair computed by Algorithm~\ref{bcd} is given
by $A = 0.15 ee^T$ and $E = I_n + 0.05 ee^T$ with error $\frac{1}{2}$. As for the Grcar matrix, allowing $E$ to be different from the identity matrix allows to reduce the error in approximating $(I_n,A)$ significantly (by a factor of two).


\small

\bibliographystyle{siam}
\bibliography{GilKS18b}

\end{document}